\documentclass{article}[12pt]

\usepackage{color}
\usepackage{amssymb}
\usepackage{mathrsfs}
\usepackage{amsmath}
\usepackage{amsfonts, amsthm, amssymb}
\usepackage{amsfonts}
\usepackage{graphicx}
\usepackage{subfigure}
\usepackage{cases}
\usepackage{booktabs}
\usepackage{float}

\newtheorem{theorem}{Theorem}[section]
\newtheorem{lemma}[theorem]{Lemma}

\newtheorem{corollary}[theorem]{Corollary}
\newtheorem{proposition}[theorem]{Proposition}
\newtheorem{problem}[theorem]{Problem}

\numberwithin{equation}{section}

\allowdisplaybreaks

\begin{document}

\title{Extremal results on degree powers \\ in some classes of graphs}

	\author{Yufei Chang, Xiaodan Chen\footnote{Corresponding author: x.d.chen@live.cn} \footnote{College of Mathematics and Information Science, Guangxi University, Nanning 530004, Guangxi, P.R. China} , Shuting Zhang}
	
 \maketitle
 
\begin{abstract}
  Let $G$ be a simple graph of order $n$ with degree sequence $(d_1,d_2,\cdots,d_n)$.
  For an integer $p>1$, let $e_p(G)=\sum_{i=1}^n d^{p}_i$ and let
  $ex_p(n,H)$ be the maximum value of $e_p(G)$ among all graphs with $n$ vertices that do not contain $H$ as a subgraph (known as $H$-free graphs).
  Caro and Yuster proposed the problem of determining the exact value of $ex_2(n,C_4)$, where $C_4$ is the cycle of length $4$.
  In this paper, we show that if $G$ is a $C_4$-free graph having $n\geq 4$ vertices and $m\leq \lfloor 3(n-1)/2\rfloor$ edges and no isolated vertices,
  then $e_p(G)\leq e_p(F_n)$, with equality if and only if $G$ is the friendship graph $F_n$.
  This yields that for $n\geq 4$, $ex_p(n,\mathcal{C}^*)=e_p(F_n)$ and $F_n$ is the unique extremal graph,
  which is an improved complement of Caro and Yuster's result on $ex_p(n,\mathcal{C}^*)$,
  where $\mathcal{C}^*$ denotes the family of cycles of even lengths.
  We also determine the maximum value of $e_p(\cdot)$
  among all minimally $t$-(edge)-connected graphs with small $t$ or among all $k$-degenerate graphs,
  and characterize the corresponding extremal graphs.
  A key tool in our approach is majorization.
\end{abstract}


\section{Introduction}
	
We only consider finite, undirected and simple graphs in this paper.
Let $G$ be a graph with vertex set $V(G)$ and edge set $E(G)$.
For a vertex $v\in V(G)$, we write $N_G(v)$ for the neighborhood of $v$ and $d_G(v)$ for the degree of $v$.
Let $\delta(G)$ denote the minimum vertex degree of $G$. 
It is well-known that $\sum_{v\in V(G)}d_G(v)=2|E(G)|$.
For an arbitrary integer $p>1$, define the {\em degree power} of $G$ as
$$e_p(G):=\displaystyle\sum_{v\in V(G)}d^{p}_G(v).$$
Hu et al. \cite{Hu} considered $e_p(G)$ for any real number $p$,
which is also known as the {\em general zeroth-order Randi\'{c} index} in the context of chemical graph theory.
Various bounds and extremal results on $e_p(G)$ of graphs without restricted structures have been established, especially for the case $p=2$;
for details one can see, e.g., \cite{Ahlswede,Caen,Chidbachian,Cioaba,Das,Fiol,Nikiforov1,Olpp,Peled,Szekely} and a comprehensive survey by Ali et al. \cite{Ali}.
	
On the other hand, along the spirit of Tur\'{a}n Theory,
Caro and Yuster \cite{Caro} introduced the problem of determining the value of $ex_p(n,H)$ for a fixed graph $H$
and characterizing the corresponding extremal graphs,
where $ex_p(n,H)$ is the maximum value of $e_p(G)$ among all graphs with $n$ vertices that do not contain $H$ as a subgraph
(also known as the $H$-free graphs).
This problem has since attracted much attention of graph theorists and many pretty results have been derived one after another;
see \cite{Bollobas2,Bollobas3,Caro} for $H$ being the complete graph $K_{r}$, \cite{Caro,Gu,Nikiforov2} for $H$ being cycles,
and \cite{Caro,Lan1,Wang} for $H$ being acyclic graphs.
For some other results we refer to \cite{Brandt,Caro,Pikhurko,Zhang1}.
Note that most of the results mentioned here just hold for sufficiently large $n$.

Let $S_n$ (resp., $C_n$) denote the star (resp., the cycle) on $n$ vertices and
$F_n$ the graph obtained from $S_n$ by adding a maximum matching on the set of its leaves,
which is named the {\it friendship graph}.
It is easy to see that $F_n$ has $\lfloor 3(n-1)/2\rfloor$ edges and no even cycles,
and that $e_p(F_n)=(n-1)^p+(n-1)\cdot 2^p$ when $n$ is odd and $e_p(F_n)=(n-1)^p+(n-2)\cdot 2^p+1$ when $n$ is even.
In \cite{Caro}, Caro and Yuster posed the following problem, which is one of the motivations for the present work.

\begin{problem}[see \cite{Caro}, Problem 6.3]\label{prob1}
Determine $ex_2(n,C_4)$. In particular, is it true that for infinitely many $n$,
$ex_2(n,C_4)=e_2(F_n)$?
\end{problem}

We now give the first theorem of this paper, which is a step towards the complete solution to Problem \ref{prob1} for $ex_p(n,C_4)$.
	
\begin{theorem}\label{th1}
Let $G$ be a $C_4$-free graph with $n\geqslant 4$ vertices and $m$ edges and $\delta(G)\geq 1$.
If $m\leq \lfloor 3(n-1)/2\rfloor$, then $e_p(G)\leq e_p(F_n)$,
with equality if and only if $G\cong F_{n}$.
\end{theorem}

Let $\mathcal{C}^{*}$ denote the family of even cycles.
Caro and Yuster \cite{Caro} proved that for sufficiently large $n$, $ex_2(n,\mathcal{C}^{*})=e_2(F_n)$ and $F_n$ is the unique extremal graph.
They further remarked that the same result holds for $ex_p(n,\mathcal{C}^{*})$ with $p>2$,
but the proof was omitted because it was rather technical.
Note that if $G$ has order $n$ and no even cycles and $e_p(G)=ex_p(n,\mathcal{C}^{*})$, then $\delta(G)\geq 1$
(otherwise, adding an edge between an isolated vertex and other arbitrary vertex of $G$ yields a new graph,
which still has no even cycles but increases $e_p(\cdot)$ strictly, a contradiction),
and $|E(G)|\leq \lfloor 3(n-1)/2\rfloor$
(since any graph on $n$ vertices with more than $\lfloor 3(n-1)/2\rfloor$ edges must contain an even cycle \cite{Caro}).
By combining these facts and Theorem \ref{th1},
we here can obtain the following corollary directly,
which is an improved complement of Caro and Yuster's result concerning $ex_p(n,\mathcal{C}^{*})$.
	
\begin{corollary}\label{coro1}
For any integer $n\geq 4$, $ex_p(n,\mathcal{C}^*)=e_p(F_{n})$ and $F_n$ is the unique extremal graph.
\end{corollary}

The problem of determining the maximum value of $e_p(G)$ taken over a class of graphs with certain special structure,
other than the $H$-free graphs, has been also well-studied in graph theory.
Many nice results of this research line have been derived for various classes of graphs
including planar graphs (and its extensions) \cite{Cook,Czap,Harant,Truszczyaski,Xu,Zhang3},
bipartite graphs \cite{Cheng,Zhang2}, and graphs with $k$ cut edges \cite{Li}.
In this paper, we consider the minimally $t$-(edge)-connected graphs and $k$-degenerate graphs.

For an integer $t\geq 1$, we say that a graph is $t$-connected (resp., $t$-edge-connected)
if removing fewer than $t$ vertices (resp., $t$ edges) always leaves the remaining graph connected,
and is minimally $t$-connected (resp., minimally $t$-edge-connected) if it is $t$-connected (resp., $t$-edge-connected)
and deleting any arbitrary chosen edge always leaves a graph which is not $t$-connected (resp., $t$-edge-connected).
Clearly, a graph is minimally 1-(edge)-connected if and only if it is a tree.
Moreover, let $K_{t,\,n-t}$ denote the complete bipartite graph with partitions of sizes $t$ and $n-t$,
which is a minimally $t$-(edge)-connected graph.
In \cite{Chen}, the authors determined the maximum value of $e_2(G)$ among all $2$-(edge)-connected graphs,
and characterized the corresponding extremal graphs.
Our next theorem shows that the same results also hold for $e_p(G)$ with an arbitrary integer $p>1$.
However, it is noted that the approach used in \cite{Chen} does not work here any more.
	
\begin{theorem}\label{th2}
(i) If $G$ is a minimally 2-connected graph with $n\geq 4$ vertices,
then $e_p(G)\leq 2\cdot(n-2)^{p}+(n-2)\cdot 2^{p}$,
with equality if and only if $G\cong K_{2,\,n-2}$.

(ii) If $G$ is a minimally 2-edge-connected graph with $n\geq 4$ vertices,
then $e_p(G)\leq \max\{e_p(F_n),e_p(K_{2,\,n-2})\}$,
with equality if and only if $G\cong F_n$ or $K_{2,\,n-2}$,
where $|V(F_n)|$ must be odd.
\end{theorem}
We remark that for fixed $p\geq 5$, if $n$ is odd and $n\geq 2p-1$, then $e_p(F_n)<e_p(K_{2,\,n-2})$ (see Appendix A (i) for its proof),
while for $p=2$ or 3 (resp., $p=4$), the necessary restriction on $n$ to guarantee this inequality is $n\geq 7$ (resp., $n\geq 9$) (verified directly by computer).

Let $W_n$ denote the graph obtained by joining each of the vertices in $C_{n-1}$ to an extra vertex,
which is named the {\it wheel graph}.
It is not hard to check that $W_n$ is a minimally 3-connected graph.
We next prove that $W_n$ and $K_{3,\,n-3}$ are the only two extremal graphs maximizing $e_p(G)$ among all minimally $3$-connected graphs.

\begin{theorem}\label{th3}
If $G$ is a minimally 3-connected graph with $n\geq 8$ vertices, then $e_p(G)\leq \max\{e_p(W_n), e_p(K_{3,\, n-3})\}$,
with equality if and only if $G\cong W_n$ or $K_{3,\,n-3}$.
\end{theorem}

We remark that for fixed $p\geq 12$, if $n\geq 2p$, then $e_p(W_n)<e_p(K_{3,\,n-3})$ (see Appendix A (ii) for its proof),
while for $2\leq p\leq 11$, the necessary restrictions on $n$ to guarantee this inequality are as follows (verified directly by computer):
\vskip2.5mm
\begin{center}
\begin{tabular}{|c|c|c|c|c|}
  \hline
  $p=2$ & $p=3$  & $p=4$  & $p=5$  & $p=6$ \\
  \hline
  $n\geq 8$ & $n\geq 9$  & $n\geq 10$  & $n\geq 12$  & $n\geq 13$ \\
  \hline
  \hline
  $p=7$ & $p=8$ & $p=9$ & $p=10$ & $p=11$ \\
  \hline
  $n\geq 15$ & $n\geq 17$ & $n\geq 19$ & $n\geq 21$ & $n\geq 23$ \\
  \hline
\end{tabular}
\end{center}

For an integer $k\geq 1$, we say that a graph is $k$-degenerate if the minimum degree in each of its subgraphs does not exceed $k$.
Clearly, all trees are $1$-degenerate.
Moreover, let $S_{n,\,k}$ denote the graph obtained from the complete graph $K_k$ on $k$ vertices
by joining each of its vertices to each of $n-k$ isolated vertices, which is apparently $k$-degenerate.
The next theorem further shows that $S_{n,\,k}$ is the unique extremal graph maximizing $e_p(G)$ among all $k$-degenerate graphs.

\begin{theorem}\label{th4}
If $G$ is a $k$-degenerate graph with $n\geqslant k+1$ vertices,
then $e_p(G)\leq k\cdot(n-1)^{p}+(n-k)\cdot k^{p}$,
with equality if and only if $G\cong S_{n,\,k}$.
\end{theorem}

It should be mentioned that a key tool for proving the above theorems is {\em majorization}.
For any two non-increasing real $n$-tuples $x=(x_1,x_2,\dots,x_n)$ and $y=$ $(y_1,y_2,\dots,y_n)$,
we say that $x$ is weakly majorized by $y$ and denote it by $x \prec_{w} y$, provided that
$$\sum_{i=1}^{s}x_i\leq\sum_{i=1}^{s}y_i,\,\, \textrm{for any integer $s$ with}\,\, 1\leq s\leq n.$$
In particular, if $\sum_{i=1}^{n}x_i=\sum_{i=1}^{n}y_i$, then we say that $x$ is majorized by $y$ and denote it by $x \prec y$.
A nice property of majorization, due to Lin et al. \cite{Lin}, is as follows, which will be used frequently in our proofs later.	

\begin{proposition}[\cite{Lin}]\label{prop1}
If $x\prec_{w} y$, then $\Vert x \Vert^{p}_{p}\leq \Vert y \Vert^{p}_{p}$,
with equality if and only if $x=y$, where $p>1$.
\end{proposition}
The rest of this paper is organized as follows.
In Section 2 we consider $C_4$-free graphs and prove Theorem \ref{th1}.
In Section 3 we consider minimally $t$-(edge)-connected graphs and prove Theorems \ref{th2} and \ref{th3}.
In Section 4 we consider $k$-degenerate graphs and prove Theorem \ref{th4}.
Some concluding remarks and open problems will be presented in the final section.

\section{$C_4$-free graphs}

In this section we will present a proof for Theorem \ref{th1} by using Proposition \ref{prop1}.
To this end, we also need the following two lemmas (see Appendix B for their proofs).

\begin{lemma}\label{lm1}
For an odd integer $n\geq 7$, let the $n$-tuples be
\begin{eqnarray}
&&\lambda_1=(n-1,2,\cdots,2,2),\nonumber\\
&&\lambda_2=\bigg(\frac{n+1}{2},\frac{n+1}{2},\frac{n+1}{2},\frac{n-1}{2},1,\ldots,1\bigg),\nonumber\\
&&\lambda_3=\big(n-q,\underbrace{q+1,\ldots,q+1}_{n-r-2},q+1-\varepsilon,\underbrace{1,\ldots,1}_r\big),\nonumber
\end{eqnarray}
where $2\leq q<\frac{n-1}{2}$, $r=\lfloor \frac{(q-1)(n-2)}{q}\rfloor$, and $\varepsilon=(q-1)(n-2)-qr$.
Then for an integer $p>1$, it follows that
(i) $\Vert\lambda_2\Vert_p^p\leq \Vert\lambda_1\Vert_p^p$ and (ii) $\Vert\lambda_3\Vert_p^p<\Vert\lambda_1\Vert_p^p$.
\end{lemma}

\begin{lemma}\label{lm12}
For an even integer $n\geq 6$, let the $n$-tuples be
\begin{eqnarray}
&&\mu_1=(n-1,2,\cdots,2,1),\nonumber\\
&&\mu_2=\bigg(\frac{n}{2}+1,\frac{n}{2},\frac{n}{2},\frac{n}{2}-1,1,\ldots,1\bigg),\nonumber\\
&&\mu_3=\big(n-q,\underbrace{q+1,\ldots,q+1}_{n-r-2},q+1-\varepsilon,\underbrace{1,\ldots,1}_r\big),\nonumber
\end{eqnarray}
where $2\leq q<\frac{n}{2}-1$, $r=\lfloor \frac{(q-1)(n-2)+1}{q}\rfloor$, and $\varepsilon=(q-1)(n-2)+1-qr$.
Then for an integer $p>1$, it follows that
(i) $\Vert\mu_2\Vert_p^p<\Vert\mu_1\Vert_p^p$ and (ii) $\Vert\mu_3\Vert_p^p<\Vert\mu_1\Vert_p^p$.
\end{lemma}

We are now ready to give a proof for Theorem \ref{th1}.\vskip 2mm

\noindent\textbf{Proof of Theorem \ref{th1}.}
Suppose that $G$ is a $C_4$-free graph with $n$ vertices and $m$ edges and $\delta(G)\geq 1$, where $n\geq 4$ and $m\leq \lfloor 3(n-1)/2\rfloor$.
Let the vertex set of $G$ be $V(G):=\{v_1,v_2,\dots,v_n\}$ and the degree sequence of $G$ be $\pi(G):=(d_1,d_2,\dots,d_n)$ with $d_i=d_G(v_i)$ for $i=1,2,\dots,n$.
Without loss of generality, we would assume that $d_1\geq d_2\geq \cdots\geq d_n$.
Note also that $F_n$ is a $C_4$-free graph with $n$ vertices and $\lfloor 3(n-1)/2\rfloor$ edges and $\delta(F_n)\geq 1$.

Clearly, to prove Theorem \ref{th1}, it suffices to show that if $G\ncong F_n$, then $e_p(G)<e_p(F_n)$.
This can be verified directly for $n=4$ or $5$.
Hence, in the following we just consider the case of $n\geq 6$.
We now assume that $G\ncong F_n$ and $e_p(G)$ is as large as possible.
This implies immediately that $d_1\leq n-2$. 
Furthermore, since $G$ is $C_4$-free, we get $|N_G(v_1)\cap N_G(v_2)|\leq 1$ and hence,
\begin{eqnarray}
d_1+d_2&=&|N_G(v_1)|+|N_G(v_2)|\nonumber\\
&=&|N_G(v_1)\cup N_G(v_2)|+|N_G(v_1)\cap N_G(v_2)|\leq n+1.\label{eq1}
\end{eqnarray}
Also, since $m\leq \lfloor 3(n-1)/2\rfloor$ and $\delta(G)\geq 1$, we obtain
$$1\leq d_n\leq \frac{2m}{n}\leq \frac{2\lfloor 3(n-1)/2\rfloor}{n}<3,$$
and thus, $d_n=1$ or 2.
We next consider the following two cases: $n$ is odd or $n$ is even.
		
{\bf Case 1:} $n$ is odd.

In this case, we have $\pi(F_n)=(n-1,2,\dots,2,2)$ and
$$\sum_{i=1}^nd_i=2m\leq 2\lfloor 3(n-1)/2\rfloor=3(n-1)=\sum \pi(F_n).$$

If $d_n=2$, then we can easily see that $\sum_{i=1}^{n-s} d_{n-i+1}\geq (n-s)\cdot 2$ holds for $1\leq s\leq n-1$,
which implies that
$$\sum_{i=1}^s d_i=2m-\sum_{i=1}^{n-s} d_{n-i+1}\leq 3(n-1)-2(n-s)=n-1+(s-1)\cdot 2.$$
This, as well as the fact that $d_1<n-1$, yields that $\pi(G)\prec_{w}\pi(F_n)$ and $\pi(G)\neq \pi(F_n)$.
Now, by Proposition \ref{prop1}, we obtain
$$e_p(G)=\Vert\pi(G)\Vert_p^p<\Vert\pi(F_n)\Vert_p^p=e_p(F_n).$$

Suppose now that $d_n=1$.
If $d_1\leq \frac{n+1}{2}$, then we consider the following $n$-tuple:
$$\pi_1:
=\bigg(\frac{n+1}{2},\frac{n+1}{2},\frac{n+1}{2},\frac{n-1}{2},1,\ldots,1\bigg).$$
It is easy to see that $\sum_{i=1}^nd_i\leq3(n-1)=\sum\pi_1$, and that for $1\leq s\leq 3$, $\sum_{i=1}^sd_i\leq s\cdot\frac{n+1}{2}$.
Also, for $4\leq s\leq n-1$, we have $\sum_{i=1}^{n-s}d_{n-i+1}\geq (n-s)$ and hence,
$$\sum_{i=1}^s d_i\leq 3(n-1)-(n-s)=3\cdot\frac{n+1}{2}+\frac{n-1}{2}+(s-4)\cdot 1.$$
We thus conclude that $\pi(G)\prec_{w}\pi_1$.
We further claim that $\pi(G)\neq \pi_1$.
Indeed, if $\pi(G)=\pi_1$, then $d_1+d_2=n+1$ and from (\ref{eq1}) we may conclude that
$N_G(v_1)\cup N_G(v_2)=V(G)$ and $N_G(v_1)\cap N_G(v_2)=\{w\}$ (say), implying that $v_1v_2\in E(G)$.
Since $G$ has no $C_4$, we see that all vertices, except $v_1$ and $v_2$, have degree at most 2,
contradicting the fact that $d_3=\frac{n+1}{2}\geq 4$ (as $n\geq 7$).
Consequently, by Proposition \ref{prop1} and Lemma \ref{lm1} (i), we get
$$e_p(G)=\Vert\pi(G)\Vert_p^p<\Vert\pi_1\Vert_p^p\leq \Vert\pi(F_n)\Vert_p^p=e_p(F_n).$$

If $d_1>\frac{n+1}{2}$, then we set $d_1=n-q$, where $q$ is an integer with $2\leq q<\frac{n-1}{2}$.
Consider the following $n$-tuple:
$$\pi_2:
=\big(n-q,\underbrace{q+1,\ldots,q+1}_{n-r-2},q+1-\varepsilon,\underbrace{1,\ldots,1}_r\big),$$
where $r=\lfloor\frac{(q-1)(n-2)}{q}\rfloor$ and $\varepsilon=(q-1)(n-2)-qr$.
Clearly, $\sum_{i=1}^nd_i\leq 3(n-1)=\sum\pi_2$.
Also, from (\ref{eq1}) it follows that $d_2\leq q+1$ and thus,
for $1\leq s\leq n-r-1$, $\sum_{i=1}^{s}d_i\leq (n-q)+(s-1)(q+1)$.
For $n-r\leq s\leq n-1$, we have $\sum_{i=1}^{n-s}d_{n-i+1}\geq (n-s)$ and thus,
$$\sum_{i=1}^s d_i\leq 3(n-1)-(n-s)=(n-q)+(n-r-2)(q+1)+(q+1-\varepsilon)+s-(n-r).$$
Now we derive that $\pi(G)\prec_{w}\pi_2$, which, as well as Proposition \ref{prop1} and Lemma \ref{lm1} (ii), would yield that
$$e_p(G)=\Vert\pi(G)\Vert_p^p\leq\Vert\pi_2\Vert_p^p<\Vert\pi(F_n)\Vert_p^p=e_p(F_n),$$
as required.
		
{\bf Case 2:} $n$ is even.

In this case, we have $\pi(F_n)=(n-1,2,\dots,2,1)$ and
$$\sum_{i=1}^nd_i=2m\leq 2\lfloor 3(n-1)/2\rfloor=3n-4=\sum \pi(F_n).$$

If $d_n=2$, then we can easily see that $\sum_{i=1}^{n-s} d_{n-i+1}\geq (n-s)\cdot 2$ holds for $1\leq s\leq n-1$,
which implies that
$$\sum_{i=1}^s d_i=2m-\sum_{i=1}^{n-s} d_{n-i+1}\leq 3n-4-2(n-s)<n-1+(s-1)\cdot 2.$$
This yields that $\pi(G)\prec_{w}\pi(F_n)$ and $\pi(G)\neq \pi(F_n)$.
Now, by Proposition \ref{prop1}, we obtain
$$e_p(G)=\Vert\pi(G)\Vert_p^p<\Vert\pi(F_n)\Vert_p^p=e_p(F_n).$$

Suppose now that $d_n=1$.
If $d_1\leq \frac{n}{2}+1$, then we consider the following $n$-tuple:
$$\pi_3:
=\bigg(\frac{n}{2}+1,\frac{n}{2},\frac{n}{2},\frac{n}{2}-1,1,\ldots,1\bigg).$$
Clearly, $\sum_{i=1}^nd_i\leq 3n-4=\sum\pi_3$.
Also, from (\ref{eq1}) it follows that $d_3\leq d_2\leq \frac{n}{2}$ and thus,
for $1\leq s\leq 3$, $\sum_{i=1}^sd_i\leq (\frac{n}{2}+1)+(s-1)\cdot\frac{n}{2}$.
Moreover, for $4\leq s\leq n-1$, we have $\sum_{i=1}^{n-s}d_{n-i+1}\geq (n-s)$ and hence,
$$\sum_{i=1}^s d_i\leq 3n-4-(n-s)=\bigg(\frac{n}{2}+1\bigg)+2\cdot\frac{n}{2}+\bigg(\frac{n}{2}-1\bigg)+(s-4)\cdot 1.$$
We thus conclude that $\pi(G)\prec_{w}\pi_3$.
Now, by Proposition \ref{prop1} and Lemma \ref{lm12} (i), we get
$$e_p(G)=\Vert\pi(G)\Vert_p^p\leq \Vert\pi_3\Vert_p^p<\Vert\pi(F_n)\Vert_p^p=e_p(F_n).$$

If $d_1>\frac{n}{2}+1$, then we set $d_1=n-q$, where $q$ is an integer with $2\leq q<\frac{n}{2}-1$.
Consider the following $n$-tuple:
$$\pi_4:
=\big(n-q,\underbrace{q+1,\ldots,q+1}_{n-r-2},q+1-\varepsilon,\underbrace{1,\ldots,1}_r\big),$$
where $r=\lfloor\frac{(q-1)(n-2)+1}{q}\rfloor$ and $\varepsilon=(q-1)(n-2)+1-qr$.
Clearly, $\sum_{i=1}^nd_i \leq3n-4=\sum\pi_4$.
Now, by the same arguments as in Case 1 one can derive that $\pi(G)\prec_{w}\pi_4$,
which, together with Proposition \ref{prop1} and Lemma \ref{lm12} (ii), would yield that
$$e_p(G)=\Vert\pi(G)\Vert_p^p\leq\Vert\pi_4\Vert_p^p<\Vert\pi(F_n)\Vert_p^p=e_p(F_n),$$
as required. This completes the proof of Theorem \ref{th1}.
$\square$

\section{Minimally $t$-(edge)-connected graphs}

In this section we will give proofs for Theorems \ref{th2} and \ref{th3}.
To this aim, we need some known results concerning the properties of minimally $t$-(edge)-connected graphs.

\begin{lemma}[\cite{Bollobas1}]\label{lm3}
Let $G$ be a minimally $t$-connected graph of order $n$.
If $n\geqslant 3t-2$ then
\begin{eqnarray}
\lvert E(G)\rvert\leqslant t(n-t). \label{eq6}
\end{eqnarray}
Furthermore, if $n\geq 3t-1$, equality holds in (\ref{eq6}) if and only if $G\cong K_{t,\,n-t}$.
\end{lemma}

\begin{lemma}[\cite{Bollobas1}]\label{lm2}
Let $G$ be a minimally $t$-connected graph. Then $\delta(G)=t$.
\end{lemma}

\begin{lemma}[\cite{Bollobas1}]\label{lm4}
A minimally 2-connected graph with more than three vertices contains no triangles.
\end{lemma}
	
\begin{lemma}[\cite{Bollobas1}]\label{lm5}
Every cycle in a minimally 3-connected graph contains at least two vertices of degree 3.
\end{lemma}

\begin{lemma}[\cite{Bondy}]\label{lm6}
Let $G$ be a minimally 2-edge-connected graph. Then $\delta(G)=2$.
\end{lemma}

A cycle $C$ of a graph $G$ is said to have a chord if there is an edge of $G$
that joins a pair of non-adjacent vertices from $C$.

\begin{lemma}[\cite{Fan}]\label{lm9}
Let $G$ be a minimally 2-edge-connected graph. Then no cycle of $G$ has a chord.
\end{lemma}

\begin{lemma}\label{lm10}
Let $G$ be a minimally 2-edge-connected graph of order $n\geq6$.
Then $|E(G)|\leq 2(n-2)$, with equality if and only if $G\cong K_{2,\,n-2}$.
\end{lemma}

\noindent\textbf{Proof.}
It suffices to prove that if $G\ncong K_{2,n-2}$, then $|E(G)|<|E(K_{2,n-2})|=2(n-2)$.
Indeed, if $G$ is 2-connected, then $G$ must be minimally 2-connected (since $G$ is minimally 2-edge-connected) and hence,
by Lemma \ref{lm3}, $|E(G)|<|E(K_{2,\,n-2})|$, as desired.
	
If $G$ is not 2-connected, then $G$ consists of some blocks, say $B_{1},B_{2},\ldots,B_{q}$ ($q\geq 2$), which connect via some cut vertices.
For $1\leq i\leq q$, since $G$ is minimally 2-edge-connected, so is $B_i$,
which yields that $|V(B_{i})|\geq 3$ (because $G$ is simple);
moreover, as above, $B_i$ must be minimally 2-connected.
For convenience, we let $n_i:=|V(B_i)|$ and $m_i:=|E(B_i)|$ and,
without loss of generality, assume that $n_1\geq \cdots \geq n_s\geq 4$ and $n_{s+1}=\cdots=n_q=3$ (here $s\geq 0$).
It is easy to see that $n=\sum_{i=1}^{s}n_i+3(q-s)-(q-1)$ and $m_{s+1}=\cdots=m_q=3$.
Furthermore, by Lemma \ref{lm3}, we have $m_i\leq 2(n_i-2)$ for $1\leq i\leq s$.
Thus, we obtain
\begin{eqnarray*}
	|E(G)|&=&\sum_{i=1}^sm_i+\sum_{i=s+1}^qm_i\\
          &\leq& \sum_{i=1}^s2(n_i-2)+3(q-s)=2n-q-s-2<2(n-2),\\
\end{eqnarray*}
where the last inequality follows from the fact that $q\geq 3$ when $s=0$ (since $n\geq 6$).
This completes the proof of Lemma \ref{lm10}.
$\square$\\

We are now ready to give the proofs for Theorems \ref{th2} and \ref{th3}.\vskip2mm

\noindent\textbf{Proof of Theorem \ref{th2}.}
We first prove (i).
Let $G$ be a minimally 2-connected graph with $V(G)=\{v_1,v_2,\ldots,v_n\}$
and let the degree sequence of $G$ be $\pi(G):=(d_1,d_2,\dots,d_n)$ with $d_i=d_G(v_i)$ for $i=1,2,\dots,n$.
Without loss of generality, assume that $d_1\geq d_2\geq \cdots\geq d_n$.
On the other hand, we see that
$$\pi(K_{2,\,n-2})=(n-2,n-2,2,\cdots,2)\,\,\, \textrm{and thus,}\,\,\, e_p(K_{2,\,n-2})=2(n-2)^p+(n-2)2^p.$$
It now suffices to show that if $G\ncong K_{2,\,n-2}$, then $e_p(G)<e_p(K_{2,\,n-2})$,
which can be verified directly for $n=4$.
So, in the following we may assume that $n\geq 5$. 
Then by Lemma \ref{lm3}, we have
$$\sum_{i=1}^nd_i=2|E(G)|<4(n-2)=\sum\pi(K_{2,\,n-2}).$$
Moreover, from Lemmas \ref{lm2} and \ref{lm4},
it follows that $d_1\leq n-2$ and thus, for $1\leq s\leq 2$, $\sum_{i=1}^sd_i\leq 2(n-2)$.
Also, since $\delta(G)=2$, for $3\leq s\leq n-1$, we have
$\sum_{i=1}^{n-s} d_{n-i+1}\geq 2(n-s)$ and hence,
$$\sum_{i=1}^s d_i=\sum_{i=1}^n d_i-\sum_{i=1}^{n-s} d_{n-i+1}<4(n-2)-2(n-s)=2\cdot(n-2)+(s-2)\cdot 2.$$
We thus conclude that $\pi(G)\prec_{w} \pi(K_{2,\,n-2})$ and $\pi(G)\neq \pi(K_{2,\,n-2})$.
Now, by Proposition \ref{prop1}, we obtain
$$e_p(G)=\Vert\pi(G)\Vert_p^p<\Vert\pi(K_{2,\,n-2})\Vert_p^p=e_p(K_{2,\,n-2}),$$
as desired.

We next prove (ii).
Let $G$ be a minimally 2-edge-connected graph with $V(G)=\{v_1,v_2,$ $\ldots,v_n\}$
and $\pi(G)=(d_1,d_2,\dots,d_n)$ with $d_i=d_G(v_i)$ for $i=1,2,\dots,n$.
Without loss of generality, assume that $d_1\geq d_2\geq \cdots\geq d_n$.
To complete the proof, we just need to prove that if $G\ncong F_n$ ($n$ is odd) and $G\ncong K_{2,\,n-2}$, then $e_p(G)<e_p(K_{2,\,n-2})$,
which can be verified directly for $n=4$ or $5$. So, in the following we may assume that $n\geq 6$.
Then by Lemma \ref{lm10}, we have
$$\sum_{i=1}^nd_i=2|E(G)|<4(n-2)=\sum\pi(K_{2,\,n-2}).$$
Moreover, we claim that $d_1\leq n-2$.
Otherwise, if $d_1=n-1$, then by Lemmas \ref{lm6} and \ref{lm9},
one can see that $(n-1)$ is even and $d_2=d_3=\cdots=d_n=2$, which implies that $G\cong F_n$ ($n$ is odd),
contradicting with the previous assumption. Our claim follows.
		
Now, based on the facts that $d_1\leq n-2$ and $\delta(G)=2$,
and using the same arguments as in the above proof for (i), we can obtain
$\pi(G)\prec_{w} \pi(K_{2,\,n-2})$ and $\pi(G)\neq \pi(K_{2,\,n-2})$.
Consequently, again by Proposition \ref{prop1}, we have
$$e_p(G)=\Vert\pi(G)\Vert_p^p<\Vert\pi(K_{2,\,n-2})\Vert_p^p=e_p(K_{2,\,n-2}).$$

This completes the proof of Theorem \ref{th2}.
$\square$\\

\noindent\textbf{Proof of Theorem \ref{th3}.}
Suppose that $G$ is a minimally 3-connected graph with $V(G)=\{v_1,v_2,\ldots,v_n\}$
and $\pi(G)=(d_1,d_2,\dots,d_n)$ with $d_i=d_G(v_i)$ for $i=1,2,\dots,n$.
Without loss of generality, assume that $d_1\geq d_2\geq \cdots\geq d_n$.
On the other hand, we have
$$\pi(K_{3,\,n-3})=(n-3,n-3,n-3,3,\ldots,3).$$
Clearly, in order to prove Theorem \ref{th3}, it suffices to show that if $G\ncong W_{n}$ and $G\ncong K_{3,\,n-3}$, then $e_p(G)<e_p(K_{3,\,n-3})$.

Indeed, since $n\geq 8$, by Lemma \ref{lm3}, we have
$$\sum_{i=1}^nd_i=2|E(G)|<6(n-3)=\sum\pi(K_{3,\,n-3}).$$
Moreover, we claim that $d_1\leq n-3$.
Otherwise, if $d_1=n-1$, then by Lemmas \ref{lm2} and \ref{lm5},
we can derive that $d_2=d_3=\cdots=d_n=3$, which implies that the subgraph of $G$ induced by $\{v_2,v_3,\ldots,v_n\}$ is isomorphic to $C_{n-1}$
(otherwise, $v_1$ is a cut vertex of $G$, a contradiction), that is $G\cong W_n$,
contradicting with the previous assumption.		
If $d_1=n-2$, then let $A:=N_G(v_1)$ and $u\in V(G)\setminus A$.
Again by Lemmas \ref{lm2} and \ref{lm5}, we can conclude that $d_G(v)=3$ holds for each $v\in A$.
We next consider the subgraph of $G$ induced by $A$, denoted by $G[A]$.
Clearly, for each $v\in A$, we see that $1\leq d_{G[A]}(v)\leq 2$;
furthermore, since $|N_G(u)\cap A|=d_G(u)\geq 3$, there are at least 4 vertices of degree one in $G[A]$.
This implies that there are at least two vertex-disjoint paths in $G[A]$ and hence,
$\{v_1,u\}$ is a 2-vertex cut in $G$, a contradiction. Our claim follows.
		
Now, based on the facts that $d_1\leq n-3$ and $\delta(G)=3$,
and using the same arguments as in the proof for (i) of Theorem \ref{th2}, we can obtain
$\pi(G)\prec_{w} \pi(K_{3,\,n-3})$ and $\pi(G)\neq \pi(K_{3,\,n-3})$.
Consequently, by Proposition \ref{prop1}, we have
$$e_p(G)=\Vert\pi(G)\Vert_p^p<\Vert\pi(K_{3,\,n-3})\Vert_p^p=e_p(K_{3,\,n-3}),$$
as desired. This completes the proof of Theorem \ref{th3}.
$\square$

\section{$k$-degenerate graphs}
	
In this section, we shall present a proof for Theorem \ref{th4}.
To this end, we need two known results regarding the properties of maximal $k$-degenerate graphs
(a $k$-degenerate graph $G$ is maximal, if for every edge $e$ of the complement of $G$, $G+e$ is not $k$-degenerate).

\begin{lemma}[\cite{Filakova}]\label{lm8}
	Let $G$ be a maximal $k$-degenerate graph with $n$ vertices, $n\geqslant k+1$. Then $G$ has $kn-{k+1\choose 2}$ edges.
\end{lemma}
	
\begin{lemma}[\cite{Filakova}]\label{lm7}
	Let $G$ be a maximal $k$-degenerate graph with $n$ vertices, $n\geqslant k+1$. Then the minimum degree of $G$ is equal to $k$.
\end{lemma}

We are now ready to give a proof for Theorem \ref{th4}.\vskip 3mm

\noindent\textbf{Proof of Theorem \ref{th4}.}
Let $G$ be a $k$-degenerate graph with $V(G)=\{v_1,v_2,\ldots,v_n\}$
and $\pi(G)=(d_1,d_2,\dots,d_n)$ with $d_i=d_G(v_i)$ for $i=1,2,\dots,n$.
Without loss of generality, we can assume that $d_1\geq d_2\geq \cdots\geq d_n$.
On the other hand, we have
$$\pi(S_{n,\,k})=(\underbrace{n-1,\dots,n-1}_{k},k,\dots,k)\,\,\, \textrm{and thus,}\,\,\, e_p(S_{n,\,k})=k(n-1)^p+(n-k)k^p.$$
Clearly, to prove Theorem \ref{th4}, it suffices to show that if $G\ncong S_{n,\,k}$, then $e_p(G)<e_p(S_{n,\,k})$.

We now assume that $G\ncong S_{n,\,k}$ and $e_p(G)$ is as large as possible.
This, as well as the fact that $e_p(G)<e_p(G+e)$ holds for any edge $e$ of the complement of $G$, yields that $G$ must be maximal.
Thus, by Lemma \ref{lm8}, we obtain
$$\sum_{i=1}^nd_i=2|E(G)|=2kn-(k+1)k=\sum \pi(S_{n,\,k}).$$
Moreover, it is easy to see that for $1\leq s\leq k$,
$\sum_{i=1}^{s}d_{i}\leq s\cdot(n-1)$.
Also, for $k+1\leq s\leq n-1$, since $\delta(G)=k$ (by Lemma \ref{lm7}), we have $\sum_{i=1}^{n-s}d_{n-i+1}\geq (n-s)\cdot k$ and hence,
\begin{equation*}
\sum_{i=1}^{s}d_{i}=\sum_{i=1}^{n}d_{i}-\sum_{i=1}^{n-s}d_{n-i+1}\leq 2kn-(k+1)k-(n-s)k=k\cdot(n-1)+(s-k)\cdot k.
\end{equation*}
We eventually derive that $\pi(G)\prec_{w} \pi(S_{n,\,k})$.
Furthermore, we have $\pi(G)\neq \pi(S_{n,\,k})$ (since $G\ncong S_{n,\,k}$).
Consequently, by Proposition \ref{prop1} we obtain
$$e_p(G)=\Vert\pi(G)\Vert_p^p<\Vert\pi(S_{n,\,k})\Vert_p^p=e_p(S_{n,\,k}),$$
as desired. This completes the proof of Theorem \ref{th4}.
$\square$

\section{Concluding remarks}

The problem of determining the exact value of $ex_2(n,C_4)$,
raised by Caro and Yuster more than 20 years ago, has not yet been solved so far.
In this paper, by using the method of majorization, we just make a step towards the complete solution to this problem by showing that
if $G$ is a $C_4$-free graph having $n\geq 4$ vertices and $m\leq \lfloor \frac{3(n-1)}{2}\rfloor$ edges and no isolated vertices,
then $e_2(G)\leq e_2(F_n)$, with equality if and only if $G$ is the friendship graph $F_n$.
However, our method seems not enough to solve this problem completely.

On the other hand, as a conjectured solution, Caro and Yuster asked whether $ex_2(n,C_4)$ $=e_2(F_n)$ holds for infinitely many $n$?
Here, we give an opposite answer to this question, which reads that
there are infinitely many $n$ such that $ex_2(n,C_4)>e_2(F_n)$.
To see this, let us recall a fact that for any prime power $q$, there exists an orthogonal polarity graph $PG(q)$ on $q^2+q+1$ vertices
such that $q+1$ vertices have degree $q$ and $q^2$ vertices has degree $q+1$ (for detailed see, e.g., \cite{He}).
Note that $PG(q)$ has no $C_4$. A simple calculation shows that for $n=q^2+q+1$ and $q\geq 5$,
$e_2(PG(q))-e_2(F_n)=q(q+1)(q-4)>0$, which implies that $ex_2(n,C_4)\geq e_2(PG(q))>e_2(F_n)$.
Naturally, an interesting problem arises:

\begin{problem} \label{prob2}
For any $C_4$-free graph $G$ on $q^2+q+1$ vertices with any integer $q\geq 5$, is it true that $e_2(G)\leq q^2(q+1)(q+2)$?
\end{problem}

\noindent If the statement in Problem \ref{prob2} is confirmed affirmatively, then $ex_2(n,C_4)=q^2(q+1)(q+2)$ holds for $n=q^2+q+1$ with prime power $q\geq 5$.

Somewhat surprising, when $p\geq 3$ and $q>1$, we have
$$ex_p(F_n)-ex_p(PG(q))=q\big(q+1\big)\cdot 2^p+q^2\big(q+1\big)\big([q^{p-2}-1][(q+1)^{p-1}-1]-1\big)>0,$$
which provides one more evidence for the assert that $ex_p(n,C_4)=e_p(F_n)$ holds for $n\geq n_0$.

In addition, by using the same tool, we determine the maximum value of $e_p(\cdot)$ among all minimally 2-(edge)-connected graphs or among all minimally $3$-connected graphs,
and characterize the corresponding extremal graphs.
It would be interesting and challenging to consider the same problem for minimally $t$-(edge)-connected graphs with general $t\geq 2$.
Formally, we pose the following problem.

\begin{problem}
For any minimally $t$-(edge)-connected graph $G$ with $n\geq f(p,t)$ vertices,
is it true that $e_p(G)\leq e_p(K_{t,\, n-t})$ with equality if and only if $G\cong K_{t,\,n-t}$?
\end{problem}

\section*{Acknowledgements}

This work was supported by the National Natural Science Foundation of China (No. 11861011).

\section*{Declarations of interest}

We declare that we have no conflict of interest.

\section*{Appendix A}

Here we prove the following two asserts:

(i) If $n$ is a odd number and $n\geq 2p-1\geq 9$, then $e_p(F_n)<e_p(K_{2,\,n-2})$.

(ii) If $n\geq 2p\geq 24$, then $e_p(W_n)<e_p(K_{3,\,n-3})$. \vskip2mm

\noindent\textbf{Proof.}
(i) Clearly, it suffices to show that $2(n-2)^{p}-(n-1)^{p}-2^p>0$ when $n\geq 2p-1\geq 9$.
Let $h_1(x):=2(x-2)^p-(x-1)^p-2^p$, where $x\geq 2p-1\geq 9$.
Then we just need to prove that $h_1(n)>0$ for $n\geq 2p-1\geq 9$.
To this aim, we consider the first derivative of $h_1(x)$ with respect to $x$:
$$h_1^\prime(x)=2p(x-2)^{p-1}-p(x-1)^{p-1}=p(x-2)^{p-1}\bigg[2-\bigg(\frac{x-1}{x-2}\bigg)^{p-1}\bigg].$$
Furthermore, noting the fact that $(1+\frac{1}{y})^y\leq e$ holds for $y\geq 0$, we get, for $x\geq 2p-1\geq 9$,
\begin{eqnarray*}
    \bigg(\frac{x-1}{x-2}\bigg)^{p-1}=
	\bigg[\bigg(1+\frac{1}{x-2}\bigg)^{x-2}\bigg]^{\frac{p-1}{x-2}}\leq e^{\frac{p-1}{x-2}}\leq e^{\frac{p-1}{2p-3}}\leq e^\frac{4}{7}\approx 1.77079<2,
\end{eqnarray*}
and hence, $h_1^\prime(x)>0$, which implies that, for $n\geq 2p-1\geq 9$,
$$h_1(n)\geq h_1(2p-1)=2(2p-3)^p-(2p-2)^p-2^p.$$
Now, in order to complete the proof, it suffices to show that $h_1(2p-1)>0$ when $p\geq 5$.
Indeed, this can be checked directly for $p=5$.
So, in the following we assume that $p\geq 6$. Consequently, we have
\begin{eqnarray*}
h_1(2p-1)&=&2^p\bigg[\bigg(\frac{2p-3}{2}\bigg)^p\bigg(2-\bigg(\frac{2p-2}{2p-3}\bigg)^p\bigg)-1\bigg]\\
&=&2^p\bigg[\bigg(\frac{2p-3}{2}\bigg)^p\bigg(2-\bigg[\bigg(1+\frac{1}{2p-3}\bigg)^{2p-3}\bigg]^{\frac{p}{2p-3}}\bigg)-1\bigg]\\
&\geq&2^p\bigg[\bigg(\frac{2p-3}{2}\bigg)^p\bigg(2-e^{\frac{p}{2p-3}}\bigg)-1\bigg]\\
&\geq&2^p\bigg[\bigg(\frac{9}{2}\bigg)^{6}\big(2-e^{\frac{2}{3}}\big)-1\bigg]>0,
\end{eqnarray*}	
as required, completing the proof of (i).

(ii) Similarly, it suffices to show that $3(n-3)^{p}-(n-1)^{p}-2\cdot3^p>0$.
Let $h_2(x):=3(x-3)^p-(x-1)^p-2\cdot3^p$, where $x\geq 2p\geq 24$.
Then we just need to prove that $h_2(n)>0$ for $n\geq 2p\geq 24$.
To this aim, we consider the first derivative of $h_2(x)$ with respect to $x$:
$$h_2^\prime(x)=3p(x-3)^{p-1}-p(x-1)^{p-1}=p(x-3)^{p-1}\bigg(3-\bigg(\frac{x-1}{x-3}\bigg)^{p-1}\bigg).$$
Furthermore, noting the fact that $(1+\frac{1}{y})^y\leq e$ holds for $y\geq 0$, we have, for $x\geq 2p\geq 24$,
\begin{eqnarray*}
    \bigg(\frac{x-1}{x-3}\bigg)^{p-1}=
	\bigg[\bigg(1+\frac{2}{x-3}\bigg)^{\frac{x-3}{2}}\bigg]^{\frac{2(p-1)}{x-3}}\leq e^{\frac{2(p-1)}{x-3}}\leq e^{\frac{2(p-1)}{2p-3}}\leq e^\frac{22}{21}\approx 2.85086<3,
\end{eqnarray*}
and hence, $h_2^\prime(x)>0$, which implies that, for $n\geq 2p\geq 24$,
$$h_2(n)\geq h_2(2p)=3(2p-3)^p-(2p-1)^p-2\cdot3^p.$$
Now, in order to complete the proof, it suffices to show that $h_2(2p)>0$ when $p\geq 12$.
Indeed, this can be checked directly by computer for $12\leq p\leq 16$.
Hence, in the following we may assume that $p\geq 17$.
Consequently, we obtain
\begin{eqnarray*}
h_2(2p)&=&3^p\bigg[\bigg(\frac{2p-3}{3}\bigg)^p\bigg(3-\bigg(\frac{2p-1}{2p-3}\bigg)^p\bigg)-2\bigg]\\
&=&3^p\bigg[\bigg(\frac{2p-3}{3}\bigg)^p\bigg(3-\bigg[\bigg(1+\frac{2}{2p-3}\bigg)^{\frac{2p-3}{2}}\bigg]^{\frac{2p}{2p-3}}\bigg)-2\bigg]\\
&\geq&3^p\bigg[\bigg(\frac{2p-3}{3}\bigg)^p\bigg(3-e^{\frac{2p}{2p-3}}\bigg)-2\bigg]\\
&\geq&3^p\bigg[\bigg(\frac{31}{3}\bigg)^{17}\big(3-e^{\frac{34}{31}}\big)-2\bigg]>0,
\end{eqnarray*}	
as desired, competing the proof of (ii).
$\qedsymbol$

\section*{Appendix B}
	
    Here, we present proofs for Lemmas \ref{lm1} and \ref{lm12} appeared in Section 2.\vskip2mm
	
	\noindent\textbf{Proof of Lemma \ref{lm1}.}
	We first prove $||\lambda_1||_p^p\geq ||\lambda_2||_p^p$ by using induction on $p$.
	Recall that
	\begin{eqnarray*}
		||\lambda_1||_p^p=(n-1)^p+(n-1)2^p\,\,\, \textrm{and}\,\,\,
		||\lambda_2||_p^p=3\bigg(\frac{n+1}{2}\bigg)^p+\bigg(\frac{n-1}{2}\bigg)^p+n-4.
	\end{eqnarray*}
	When $p=2$, we have
	\begin{eqnarray*}
		||\lambda_1||_2^2-||\lambda_2||_2^2=(n-1)^2+4(n-1)-\bigg(\frac{3(n+1)^2}{4}+\frac{(n-1)^2}{4}+n-4\bigg)=0,
	\end{eqnarray*}
	as desired. Suppose now that $||\lambda_1||_p^p\geq ||\lambda_2||_p^p$ holds for $p=k\geq 2$.
	We next consider the case when $p=k+1$. By using the inductive hypothesis, we get
	\begin{eqnarray*}
		||\lambda_1||_{k+1}^{k+1}
		&=&(n-1)||\lambda_1||_{k}^{k}-(n-3)(n-1)2^k\\
		&\geq&(n-1)||\lambda_2||_{k}^{k}-(n-3)(n-1)2^k\\
		&=&||\lambda_2||_{k+1}^{k+1}+\phi(n),
	\end{eqnarray*}
	where
	$$\phi(n):=\frac{3(n-3)}{2}\bigg(\frac{n+1}{2}\bigg)^k+\bigg(\frac{n-1}{2}\bigg)^{k+1}+(n-2)(n-4)-(n-3)(n-1)2^k.$$
	Furthermore, for $n\geq 7$ and $k\geq 2$, we have
	\begin{eqnarray}
		&&2n^3-19n^2+58n-53\geq 0\nonumber\\
		&\Rightarrow& \frac{3(n-3)}{2}\bigg(\frac{n+1}{4}\bigg)^2+\frac{n-1}{2}\bigg(\frac{n-1}{4}\bigg)^{2}\geq (n-3)(n-1)\nonumber\\
		&\Rightarrow& \frac{3(n-3)}{2}\bigg(\frac{n+1}{4}\bigg)^k+\frac{n-1}{2}\bigg(\frac{n-1}{4}\bigg)^{k}+\frac{(n-2)(n-4)}{2^k}\geq (n-3)(n-1)\nonumber\\
		&\Rightarrow& \frac{3(n-3)}{2}\bigg(\frac{n+1}{2}\bigg)^k+\bigg(\frac{n-1}{2}\bigg)^{k+1}+(n-2)(n-4)\geq (n-3)(n-1)2^k\nonumber\\
		&\Rightarrow& \phi(n)\geq 0.\nonumber
	\end{eqnarray}
	This proves $||\lambda_1||_{k+1}^{k+1}\geq ||\lambda_2||_{k+1}^{k+1}$ and
	by the principle of induction, we complete the proof of the first part.
	
	We next prove $||\lambda_1||_p^p>||\lambda_3||_p^p$.
	Recall also that
	\begin{eqnarray*}
		||\lambda_3||_p^p=\big(n-q\big)^p+\big(n-r-2\big)\big(q+1\big)^p+\big(q+1-\varepsilon\big)^p+r,
	\end{eqnarray*}
	where $2\leq q<\frac{n-1}{2}$, $r=\lfloor\frac{(q-1)(n-2)}{q}\rfloor$, and $\varepsilon=(q-1)(n-2)-qr$.
	Consider the following function $f_p(x)$ defined on the interval $[a,\,b]$, where $a=\frac{(q-1)(n-2)}{q}-1$, $b=\frac{(q-1)(n-2)}{q}$, and
	\begin{eqnarray*}
		f_p(x):=\big(n-q\big)^p+\big(n-x-2\big)\big(q+1\big)^p+\big(q+1+qx-(q-1)(n-2)\big)^p+x.
	\end{eqnarray*}
	It is easy to check that
	\begin{eqnarray*}
		f_p^{\prime\prime}(x)=(p-1)pq^2\big(q+1+qx-(q-1)(n-2)\big)^{p-2}\geq(p-1)pq^2>0,
	\end{eqnarray*}
	which implies that
	\begin{eqnarray}
		||\lambda_3||_p^p=f_p(r)&\leq& \max\big\{f_p(a), f_p(b)\big\}\nonumber\\
&=&(n-q)^p+\frac{n+q-2}{q}(q+1)^p+\frac{(q-1)(n-2)}{q}. \label{eqA1}
	\end{eqnarray}
	We now consider the following two cases:
	
	\textbf{Case 1.} $q=2$
	
	In this case, noting that $n$ is an odd integer with $2\leq q<\frac{n-1}{2}$, we have $n\geq 2q+3=7$.
	If $p=2$, then by (\ref{eqA1}) we obtain
	$$||\lambda_1||_{2}^{2}-||\lambda_3||_{2}^{2}=n-6>0.$$
	If $p\geq 3$, then again by (\ref{eqA1}) we get
	\begin{eqnarray}
		||\lambda_1||_{p}^{p}-||\lambda_3||_{p}^{p}
		&\geq& (n-1)^{p}-(n-2)^{p}-(3^{p}-2^{p})\frac{n}{2}+(2^{p}-1)(\frac{n}{2}-1)\nonumber\\
		&\geq& p(n-2)^{p-1}-(3^{p}-2^{p})\frac{n}{2}\nonumber\\
		&\geq& p\cdot 5^{p-1}-\frac{7}{2}(3^{p}-2^{p})\nonumber\\
		&\geq& 3\cdot 5^2-\frac{7}{2}(3^{3}-2^{3})>0,\nonumber
	\end{eqnarray}
	where the second inequality follows from a direct application of Lagrange's mean value theorem to the function $x^{p}$ on the interval $[n-2,n-1]$
	and the fact that $(2^{p}-1)(\frac{n}{2}-1)>0$,
	while the third one follows from the fact that $p(n-2)^{p-1}-(3^{p}-2^{p})n/2$ is increasing with respect to $n$ when $n\geq 7$ and $p\geq 3$.
	
	\textbf{Case 2.} $3\leq q<\frac{n-1}{2}$
	
	In this case, we have $n\geq 2q+3\geq 9$.
	If $p=2$, then by (\ref{eqA1}) we get
	\begin{eqnarray*}
		||\lambda_1||_2^2-||\lambda_3||_2^2\geq (q-1)\big[n-2(q+1)\big]>0.
	\end{eqnarray*}
	If $p\geq 3$, then again by (\ref{eqA1}), we have
	\begin{eqnarray}
		||\lambda_1||_{p}^{p}-||\lambda_3||_{p}^{p}
		&\geq& (n-1)^{p}-(n-q)^{p}-\frac{n+q-2}{q}(q+1)^p\nonumber\\
        &&+(n-1)2^{p}-\frac{(q-1)(n-2)}{q}\nonumber\\
		&\geq& p(q-1)(n-q)^{p-1}-\frac{n+q-2}{q}(q+1)^p\nonumber\\
		&\geq& \frac{1}{q}\big[pq(q-1)(q+3)^{p-1}-(3q+1)(q+1)^p\big]\nonumber\\
		&\geq& \frac{1}{q}\big[3q(q-1)(q+3)^2-(3q+1)(q+1)^3\big]\nonumber\\
		&=& \frac{1}{q}\big(5q^3-3q^2-33q-1\big)>0,\nonumber
	\end{eqnarray}
	where the third inequality follows from the fact that $p(q-1)(n-q)^{p-1}-(n+q-2)(q+1)^p/q$ is increasing with respect to $n$
	when $n\geq 2q+3$, $q\geq 3$ and $p\geq 3$.
	
	This completes the proof of Lemma \ref{lm1}.
	$\qedsymbol$\\

	\noindent\textbf{Proof of Lemma \ref{lm12}.}
	We first prove $||\mu_1||_p^p> ||\mu_2||_p^p$ by using induction on $p$.
	Recall that
	\begin{eqnarray*}
		&&||\mu_1||_p^p=(n-1)^p+(n-2)2^p+1,\nonumber\\
		&&||\mu_2||_p^p=\bigg(\frac{n}{2}+1\bigg)^p+2\bigg(\frac{n}{2}\bigg)^p+\bigg(\frac{n}{2}-1\bigg)^p+n-4.
	\end{eqnarray*}
	When $p=2$, we have
	\begin{eqnarray*}
		||\mu_1||_2^2-||\mu_2||_2^2=n-4>0,
	\end{eqnarray*}
	as desired. Suppose now that $||\mu_1||_p^p>||\mu_2||_p^p$ holds for $p=k\geq 2$.
	We next consider the case when $p=k+1$. By using the inductive hypothesis, we get
	\begin{eqnarray*}
		||\mu_1||_{k+1}^{k+1}
		&=&(n-1)||\mu_1||_{k}^{k}-(n-3)(n-2)2^k-(n-2)\\
		&>&(n-1)||\mu_2||_{k}^{k}-(n-3)(n-2)2^k-(n-2)\\
		&=&||\mu_2||_{k+1}^{k+1}+\varphi(n),
	\end{eqnarray*}
	where
	$$\varphi(n):=\frac{n-4}{2}\bigg(\frac{n}{2}+1\bigg)^k+(n-2)\bigg(\frac{n}{2}\bigg)^{k}+\frac{n}{2}\bigg(\frac{n}{2}-1\bigg)^k+(n-2)(n-5)-(n-3)(n-2)2^k.$$
	Furthermore, for $n\geq 6$ and $k\geq 2$, we have
	\begin{eqnarray}
		&&n^3-10n^2+38n-52> 0\nonumber\\
		&\Rightarrow& \frac{n-4}{2}\bigg(\frac{n+2}{4}\bigg)^2+(n-2)\bigg(\frac{n}{4}\bigg)^{2}+\frac{n}{2}\bigg(\frac{n-2}{4}\bigg)^2> (n-3)(n-2)\nonumber\\
		&\Rightarrow& \frac{n-4}{2}\bigg(\frac{n+2}{4}\bigg)^k+(n-2)\bigg(\frac{n}{4}\bigg)^{k}+\frac{n}{2}\bigg(\frac{n-2}{4}\bigg)^k+\frac{(n-2)(n-5)}{2^k}\nonumber\\
        &&> (n-3)(n-2)\nonumber\\
		&\Rightarrow& \frac{n-4}{2}\bigg(\frac{n+2}{2}\bigg)^k+(n-2)\bigg(\frac{n}{2}\bigg)^{k}+\frac{n}{2}\bigg(\frac{n-2}{2}\bigg)^k+(n-2)(n-5)\nonumber\\
        &&> (n-3)(n-2)2^k\nonumber\\
		&\Rightarrow& \varphi(n)> 0.\nonumber
	\end{eqnarray}
	This proves $||\mu_1||_{k+1}^{k+1}> ||\mu_2||_{k+1}^{k+1}$ and
	by the principle of induction, we complete the proof of the first part.
	
	We next prove $||\mu_1||_p^p>||\mu_3||_p^p$.
	Recall also that
	\begin{eqnarray*}
		||\mu_3||_p^p=\big(n-q\big)^p+\big(n-r-2\big)\big(q+1\big)^p+\big(q+1-\varepsilon\big)^p+r,
	\end{eqnarray*}
	where $2\leq q<\frac{n}{2}-1$, $r=\lfloor\frac{(q-1)(n-2)+1}{q}\rfloor$, and $\varepsilon=(q-1)(n-2)+1-qr$.
	Consider the following function $g_p(x)$ on the interval $[c,\,d]$, where $c=\frac{(q-1)(n-2)+1}{q}-1$, $d=\frac{(q-1)(n-2)+1}{q}$, and
	\begin{eqnarray*}
		g_p(x):=\big(n-q\big)^p+\big(n-x-2\big)\big(q+1\big)^p+\big(q+qx-(q-1)(n-2)\big)^p+x.
	\end{eqnarray*}
	It is easy to check that
	\begin{eqnarray*}
		g_p^{\prime\prime}(x)=(p-1)pq^2\big(q+qx-(q-1)(n-2)\big)^{p-2}\geq(p-1)pq^2>0,
	\end{eqnarray*}
	which implies that
	\begin{eqnarray}
		||\mu_3||_p^p=g_p(r)&\leq& \max\big\{g_p(c), g_p(d)\big\}\nonumber\\
        &=&(n-q)^p+\frac{n+q-3}{q}(q+1)^p+\frac{(q-1)(n-2)+1}{q}. \label{eqA2}
	\end{eqnarray}
	We now consider the following two cases:
	
	\textbf{Case 1.} $q=2$
	
	In this case, noting that $n$ is an even integer with $2\leq q<\frac{n}{2}-1$, we have $n\geq 2q+4=8$.
	If $p=2$, then by (\ref{eqA2}) we get
	\begin{eqnarray*}
		||\mu_1||_2^2-||\mu_3||_2^2\geqslant n-5>0.
	\end{eqnarray*}
	If $p\geqslant3$, then again by (\ref{eqA2}) we obtain
	\begin{eqnarray}
		||\mu_1||_{p}^{p}-||\mu_3||_{p}^{p}
		&\geq& (n-1)^{p}-(n-2)^{p}-(3^{p}-2^{p})\frac{n-1}{2}+(2^{p}-1)\frac{n-3}{2}\nonumber\\
		&\geq& p(n-2)^{p-1}-(3^{p}-2^{p})\frac{n-1}{2}\nonumber\\
		&\geq& p\cdot 6^{p-1}-\frac{7}{2}(3^{p}-2^{p})\nonumber\\
		&\geq& 3\cdot 6^2-\frac{7}{2}(3^{3}-2^{3})>0,\nonumber
	\end{eqnarray}
	where the second inequality follows from a direct application of Lagrange's mean value theorem to the function $x^{p}$ on the interval $[n-2,n-1]$
	and the fact that $(2^{p}-1)(n-3)/2>0$,
	while the third one follows from the fact that $p(n-2)^{p-1}-(3^{p}-2^{p})(n-1)/2$ is increasing with respect to $n$ when $n\geq 8$ and $p\geq 3$.
	
	\textbf{Case 2.} $3\leq q<\frac{n}{2}-1$
	
	In this case, we have $n\geq 2q+4\geq 10$.
	If $p=2$, then by (\ref{eqA2}) we get
	\begin{eqnarray*}
		||\mu_1||_2^2-||\mu_3||_2^2&\geq& (q-1)[n-(2q+1)]>0.
	\end{eqnarray*}
	If $p\geq 3$, then again by (\ref{eqA2}), we obtain
	\begin{eqnarray}
		||\mu_1||_{p}^{p}-||\mu_3||_{p}^{p}
		&\geq& (n-1)^{p}-(n-q)^{p}-\frac{n+q-3}{q}(q+1)^p\nonumber\\
        &&+(n-2)2^{p}-\frac{(q-1)(n-3)}{q}\nonumber\\
		&\geq& p(q-1)(n-q)^{p-1}-\frac{n+q-3}{q}(q+1)^p\nonumber\\
		&\geq& \frac{1}{q}\big[pq(q-1)(q+4)^{p-1}-(3q+1)(q+1)^p\big]\nonumber\\
		&\geq& \frac{1}{q}\big[3q(q-1)(q+4)^2-(3q+1)(q+1)^3\big]\nonumber\\
		&=& \frac{1}{q}\big(11q^3+12q^2-54q-1\big)>0,\nonumber
	\end{eqnarray}
	where the third inequality follows from the fact that $p(q-1)(n-q)^{p-1}-(n+q-3)(q+1)^p/q$
	is increasing with respect to $n$ when $n\geq 2q+4$, $q\geq 3$ and $p\geq 3$.
	
	This completes the proof of Lemma \ref{lm12}.
	$\qedsymbol$

\end{document}